%% file: bote-cavity.tex
\documentclass[journal,transmag]{IEEEtran}
\usepackage{cite}
\usepackage{amssymb}
\usepackage[pdftex]{graphicx}
\usepackage{amsfonts}
\usepackage{lineno}
\usepackage[cmex10]{amsmath}
\usepackage{array}
\usepackage{mathrsfs}
\usepackage{stfloats}
\usepackage{url}
\usepackage{epstopdf}
\usepackage[square,numbers,sort&compress]{natbib}
\usepackage{subfigure}
\usepackage{subeqnarray}
\usepackage{cases}
\usepackage{color}

\input latex_defs.tex

\allowdisplaybreaks
\input latex_defs.tex
\begin{document}
\title{Numerical Eigensolver for Solving Eigenmodes of Cavity Resonators Filled With both Electric and Magnetic Lossy, Anisotropic Media}
\author{\IEEEauthorblockN{Wei Jiang\IEEEauthorrefmark{1}, Jie Liu\IEEEauthorrefmark{2} and Shiling Zheng\IEEEauthorrefmark{3}}
\IEEEauthorblockA{\IEEEauthorrefmark{1} School of Mechatronics Engineering, Guizhou Minzu University, Guiyang 550025, China
}
\IEEEauthorblockA{\IEEEauthorrefmark{2} Institute of Electromagnetics
and Acoustics, Xiamen University, Xiamen 361005, China
}
\IEEEauthorblockA{\IEEEauthorrefmark{3} School of Physics Science and Information Technology, Liaocheng University, Liaocheng 252000, China
}


}

\IEEEtitleabstractindextext{
\begin{abstract}
This article presents the numerical eigensolver to find the resonant frequencies of 3-D closed cavity resonators filled with both electric and magnetic lossy, anisotropic media. By introducing a dummy variable with zero value in the 3-D linear vector Maxwell eigenvalue problem for the electric field, we enforce the divergence-free condition for electric flux density in a weak sense. In addition, by introducing a dummy variable with constant value in the 3-D linear vector Maxwell eigenvalue problem for the magnetic field, we enforce the divergence-free condition for magnetic flux density in a weak sense.
Moreover, it is theoretically proved that the novel method of introducing dummy variables can be free of all the spurious modes in solving eigenmodes of the 3-D closed cavity problem. Numerical experiments show that the numerical eigensolver supported by this article can eliminate all the spurious modes, including spurious dc modes.
\end{abstract}
}

\maketitle
\IEEEdisplaynontitleabstractindextext
\IEEEpeerreviewmaketitle

\section{Introduction}
\IEEEPARstart{I}{t} is well known that finding several resonant eigenmodes of 3-D closed cavity problems is one of the classic problems in computational electromagnetics. It is also known that spurious modes will appear in numerical results if the numerical calculation method cannot maintain the physical properties of the electromagnetic field correctly.
Now, we know that divergence-free condition supported by Gauss's law is very important to numerical simulations of the resonant cavity problem.
If the divergence-free condition is artificially ignored in the numerical method, then a lot of spurious dc modes will be introduced in the numerical results.
For example, some early articles \cite{Lee1992,Wang1996,Kanai2000} about employing edge elements to simulate 3-D closed cavity problem do not enforce the divergence-free condition, as a result, numerical methods given in these articles do not remove spurious dc modes. When an inappropriate numerical method is used to solve the resonant cavity problem, then spurious nonzero modes will be introduced in the numerical results. In order to eliminate these spurious nonzero and dc modes in the numerical results together, we must select an appropriate computational method to solve 3-D closed cavity problem.

We know that the finite element method (FEM) is an important numerical computational method for solving partial differential equations (PDEs). In recent years, Jiang \emph{et. al.} \cite{Jiang2016,Jiang2019} have applied mixed finite element method (MFEM) to solve the 3-D resonant cavity problem filled with electric and/or magnetic lossless media, excluding the case that both electric and magnetic lossy media. For all of these cavity problems, the MFEM can remove all the spurious modes, including spurious dc modes,  because in the MFEM there is a Lagrange multiplier to enforce the divergence-free condition in a weak sense.
Based upon the idea of the MFEM, Liu \emph{et. al.} \cite{Liu2017} propose a two grid method to solve the abovementioned 3-D cavity problem. Moreover, the method given in \cite{Liu2017} can remove all the spurious modes, including spurious dc modes.
For the 3-D cavity problem filled with fully conductive media in the whole cavity, Jiang \emph{et. al.} \cite{Jiang2017} successfully employ edge element of the lowest order and the standard linear element to solve this problem.
Note that the medium considered in \cite{Jiang2017} is conductive loss, not dielectric loss.

In 2020, Jiang \emph{et. al.} \cite{Jiang2020} provide three numerical eigensolvers to solve the 3-D resonant cavity problem. In \cite{Jiang2020}, the projection method can solve 3-D resonant cavity problem filled with both electric and magnetic lossy media. However, this method need to obtain the null space of a given sparse matrix.
In the projection method, a completed singular value deposition method is used to find a set of
normalized orthogonal bases of a null space. Since the cost of completed singular value decomposition is very huge, the projection method given in \cite{Jiang2020} is not efficient. For the 3-D resonant cavity problem filled with both electrical and magnetic lossy, anisotropic media, the penalty method given in \cite{Jiang2020} cannot remove all the spurious modes.
Furthermore, the augmented method given in \cite{Jiang2020} cannot guarantee that the numerical modes obtained by the augmented method are all physical.
Based on these reasons, this article continues to investigate the 3-D resonant cavity problem filled with both electrical and magnetic lossy, anisotropic media.
In this article, on one hand, by introducing a dummy variable with zero value in the 3-D vector Maxwell's eigenvalue problem for the electric field, we enforce the divergence-free condition for electric flux density in a weak sense. On other hand, by introducing a dummy variable with constant value in the 3-D vector Maxwell's eigenvalue problem for the magnetic field, we enforce the divergence-free condition for magnetic flux density in a weak sense.
Moreover, it is theoretically proved that the novel method of introducing dummy variables can be free of all the spurious modes, including spurious dc modes.

\indent The outline of the paper is as follows. The well posed governing equations for the 3-D closed cavity problem filled with both electric and magnetic lossy media are given in Section \uppercase\expandafter{\romannumeral2}. In Section \uppercase\expandafter{\romannumeral3},
we present two mixed variational forms and provide the FEM to solve them. In Section \uppercase\expandafter{\romannumeral4}, two numerical experiments are carried out to verify that FEM is free of all the spurious modes, including spurious dc modes.


\section{The Governing Equations of 3-D Resonant Cavity Problem}
Let $\Omega$ be the geometry of a given cavity resonator, which is usually a bounded domain in $\mathbb{R}^3$.
 The boundary of the cavity resonator is denoted by $\p\Omega$. Note that $\p\Omega$ may be disconnected, in this case,
 the cavity $\Omega$ is not a simply connected manifold in $\mathbb{R}^3$.
 Let $\^n$ be the unit outward normal vector on the boundary $\p\Omega$. The permittivity and permeability in vacuum are denoted by $\ep_{0}$ and $\mu_{0}$, respectively. In this article, we mainly focus on finding the resonant eigenmodes of the resonant cavity filled with both electric and magnetic lossy, anisotropic media. Hence, the relative permittivity tensor $\d{\ep}_{r}$ and the relative permeability $\d{\mu}_{r}$ of an anisotropic medium are not Hermitian \cite{Chew1990book}, that is to say that both $\d{\ep}^{\dag}\neq\d{\ep}$ and
 $\d{\mu}_{r}^{\dag}\neq\d{\mu}_{r}$ holds, where the symbol $\dag$ stands for complex conjugate transpose.

According to the classic electromagnetic theory \cite{Balanis}, we know that the governing equations of the 3-D closed cavity problem are the source-free Maxwell's equations of the first order.
In addition, the governing equations of the 3-D closed cavity problem include perfect electric conductor boundary condition. On one hand, eliminating the magnetic field $\H$ from the source-free Maxwell's equations of the first order, we can obtain a linear PDE system of the second order with respect to the electric field $\E$, and that is
\begin{subequations} \label{eq:2}
\begin{numcases}{}
  \curl\Big({\d{\mu}_{r}^{-1}}\curl\E\Big) =\Lambda{\d{\ep}}_{r}\E~~\text{in}~\Omega\label{eq:2a}\\
  \div\big({\d{\ep}}_{r}\E\big)= 0 ~~\text{in}~\Omega\label{eq:2b}\\
  \^n\times\E = {\bf{0}}~~\text{on}~\p\Omega,\label{eq:2c}
\end{numcases}
\end{subequations}
where $(\Lambda,\E)$ with $\E\neq{\bf{0}}$ are unknown.

In equation (\ref{eq:2a}), $\Lambda=\omega^2\ep_{0}\mu_{0}$ is the square of the wavenumber in vacuum. Since the material in the cavity is lossy, then the spectral point $\Lambda$ is made up of countable complex numbers. On the other hand, eliminating the electric field $\E$ from the source-free Maxwell's equations of the first order, one can obtain a linear PDE system of the second order with respect to the magnetic field $\H$, and that is
\begin{subequations} \label{eq:3}
\begin{numcases}{}
  \curl\Big({\d{\ep}_{r}^{-1}}\curl\H\Big) =\Lambda{\d{\mu}}_{r}\H~~\text{in}~\Omega\label{eq:3a}\\
  \div\big({\d{\mu}}_{r}\H\big)= 0 ~~\text{in}~\Omega\label{eq:3b}\\
  \^n\times({\d{\ep}_{r}^{-1}}\curl\H\Big) = {\bf{0}}~~\text{on}~\p\Omega\label{eq:3d}\\
   \^n\cdot({\d{\mu}}_{r}\H\big) = {{0}}~~\text{on}~\p\Omega,\label{eq:3c}
\end{numcases}
\end{subequations}
where $(\Lambda,\H)$ with $\H\neq{\bf{0}}$ are unknown. The nonzero eigenmodes between PDEs (\ref{eq:2}) and PDEs (\ref{eq:3}) are equivalent if the electromagnetic field is at least twice differentiable. When the geometry $\Omega$ of the resonant cavity has a very complex geometric topology such that both PDEs (\ref{eq:2}) and PDEs (\ref{eq:3}) have physical dc modes, in this case, it is worthwhile to point out that these dc modes are not equivalent between PDEs (\ref{eq:2}) and PDEs (\ref{eq:3}). For details, please see \cite{Jiang2019p}.

It is known that a three dimensional vector function has three scalar function components, therefore there are four single scalar equations in (\ref{eq:2a})-(\ref{eq:2b}) and (\ref{eq:3a})-(\ref{eq:3b}), respectively. That is to say that PDEs (\ref{eq:2}) and PDEs (\ref{eq:3}) are two linear overdetermined systems. Based on this fact, we introduce two new scalar functions $p$ and $q$ in (\ref{eq:2a}) and (\ref{eq:3a}) respectively, such that PDEs (\ref{eq:2}) and PDEs (\ref{eq:3}) are well posed.

The well posed eigenvalue problem associated with PDEs (\ref{eq:2}) is as follow:
\begin{subequations} \label{eqq:2}
\begin{numcases}{}
  \curl\Big({\d{\mu}_{r}^{-1}}\curl\E\Big)+\alpha\nabla p =\Lambda{\d{\ep}}_{r}\E~~\text{in}~\Omega\label{eqq:2a}\\
  \div\big({\d{\ep}}_{r}\E\big)= 0 ~~\text{in}~\Omega\label{eqq:2b}\\
  \^n\times\E = {\bf{0}}~~\text{on}~\p\Omega\label{eqq:2c}\\
  p=0~~\text{on}~\p\Omega,\label{eqq:2d}
\end{numcases}
\end{subequations}
where $(\Lambda,\E,p)$ with $\E\neq{\bf{0}}$ is unknown and $\alpha$ is a given positive constant.

The well posed eigenvalue problem associated with PDEs (\ref{eq:3}) is as follow:
\begin{subequations} \label{eqq:3}
\begin{numcases}{}
  \curl\Big({\d{\ep}_{r}^{-1}}\curl\H\Big)+\beta\nabla q =\Lambda{\d{\mu}}_{r}\H~~\text{in}~\Omega\label{eqq:3a}\\
  \div\big({\d{\mu}}_{r}\H\big)= 0 ~~\text{in}~\Omega\label{eqq:3b}\\
 \^n\times\Big({\d{\ep}_{r}^{-1}}\curl\H\Big) = {\bf{0}}~~\text{on}~\p\Omega\label{eqq:3c}\\
   \^n\cdot\big({\d{\mu}}_{r}\H\big) = {{0}}~~\text{on}~\p\Omega\label{eqq:3d}\\
   \pd{q}{n}=0~~\text{on}~\p\Omega,\label{eqq:3e}
\end{numcases}
\end{subequations}
where $\beta$ is a given positive constant. Here, our aim is to find $(\Lambda,\H,q)$ with $\H\neq{\bf{0}}$ such that PDEs (\ref{eqq:3}) holds. In the next section, we give one way to select the constant $\alpha$ in (\ref{eqq:2a}) and the constant $\beta$ in (\ref{eqq:3a}).\\
\indent Now we prove that the scalar function $p$ in PDEs (\ref{eqq:2}) is a dummy variable with zero value. As a matter of fact,
by taking the divergence at the both of (\ref{eqq:2a}) and then use (\ref{eqq:2b}) and (\ref{eqq:2d}), at last we get the following PDE
\begin{subequations} \label{harm}
\begin{numcases}{}
  \nabla^2 p=0~\text{in}~\Omega\label{harm1}\\
  p=0~~\text{on}~\p\Omega.\label{harm2}
\end{numcases}
\end{subequations}
It follows that $p$ is a harmonic function from (\ref{harm1}). According to the classic extremum property of the harmonic function \cite{evans1998}, it is easy to obtain the conclusion that the function $p$ vanishes in the whole domain $\Omega$. As a consequence, the solution $(\Lambda,\E,0)$ to (\ref{eqq:2}) is the solution $(\Lambda,\E)$ to (\ref{eq:2}). Obviously, the solution $(\Lambda,\E)$ to (\ref{eq:2}) is also the solution $(\Lambda,\E,p)$ to (\ref{eqq:2}) as long as we take $p=0~\text{in}~\Omega$. Therefore, the solution between (\ref{eq:2}) and (\ref{eqq:2})
is equivalent. Next we prove the scalar function $q$ in PDEs (\ref{eqq:3}) is a dummy variable with any constant value. As a matter of fact, by taking the divergence at the both of (\ref{eqq:3a}) and then use (\ref{eqq:3b}) and (\ref{eqq:3e}), at last we get the following PDE
\begin{subequations} \label{harmb}
\begin{numcases}{}
  \nabla^2 q=0~\text{in}~\Omega\label{harmb1}\\
  \pd{q}{n}=0~~\text{on}~\p\Omega.\label{harmb2}
\end{numcases}
\end{subequations}
Using the first Green's theorem of scalar function, one can obtain
\begin{equation}\label{thm1}
    \int_{\Omega}\grad q\cdot \grad q \textrm{d}\Omega +\int_{\Omega}q\nabla^2 q\textrm{d}\Omega=\int_{\p\Omega}q\pd{q}{n}\textrm{d}S.
\end{equation}
Substituting  (\ref{harmb1}) and (\ref{harmb2}) into (\ref{thm1}), we have
\begin{equation}\label{tht1}
\int_{\Omega}\grad q\cdot \grad q \textrm{d}\Omega=0.
\end{equation}
The above equation (\ref{tht1}) implies $\grad q={\bf{0}}$ in $\Omega$. Therefore the function $q$ is an arbitrary constant. That is to say that the solution $(\Lambda,\H,C)$ to (\ref{eqq:3}) is the solution $(\Lambda,\H)$ to (\ref{eq:3}), where $C$ is an arbitrary constant. Obviously, the solution $(\Lambda,\H)$ to (\ref{eq:3}) is also the solution $(\Lambda,\H, q)$ to (\ref{eqq:3}) as long as we take $q=C~\text{in}~\Omega$. Therefore, the solution between (\ref{eq:3}) and (\ref{eqq:3})
is also equivalent.



\section{Mixed Finite Element Method Discretization}
In this section, the MFEM in computational mathematics is applied to solve PDEs (\ref{eqq:2}) and PDEs (\ref{eqq:3}). In order to give the mixed variational forms associated with PDEs (\ref{eqq:2}) and PDEs (\ref{eqq:3}), at first, we need to introduce the following Hilbert spaces:
\begin{gather*}
     L^{2}(\Omega)=\Big\{f: \int_{\Omega}|f(x,y,z)|^2 \mathrm {d}x\mathrm{d}y\mathrm{d}z<+\infty\Big\}\\
    H^{1}(\Omega)=\big\{f\in{L^2(\Omega):\grad f\in{(L^2(\Omega))^3}}\big\}\\
    H_{0}^1(\Omega)=\big\{f\in{H^1(\Omega)}:f=0\textrm{~on~}\p\Omega\big\}\\
    \mathbb{H}(\textrm{curl},\Omega)=\big\{\F\in{(L^2(\Omega))^3}: \curl{\F}\in{(L^2(\Omega))^3}\big\}\\
    \mathbb{H}_{0}(\textrm{curl},\Omega)=\big\{\F\in{\mathbb{H}(\textrm{curl},\Omega)}: \^n\times{\F}=0\textrm{~on~}\p\Omega\big\}
\end{gather*}
Next, let us introduce the mixed variational forms associated with PDEs (\ref{eqq:2}) and PDEs (\ref{eqq:3}).
\subsection{Mixed Variational Forms}
Using Green formulas, the mixed variational form associated with PDEs (\ref{eqq:3}) reads as:
\begin{small}
	\begin{subequations} \label{eqnhtt:3}
		\begin{numcases}{}
			\textrm{Seek~} \Lambda\in{\mathbb{C}},~{\bf{0}}\neq\H\in{\mathbb{H}(\textrm{curl},\Omega)},~q\in{H^1(\Omega)},\textrm{~such that~}\nonumber\\
			a_{1}(\H,\F)+\beta b_{1}(\F,q)= \Lambda d_{1}(\H,\F),~\forall~\F\in{\mathbb{H}(\textrm{curl},\Omega)}\label{eqnhtt:3a}\\
			c_{1}(\H,\upsilon) = 0,~~\forall\upsilon\in{H^{1}(\Omega)}\label{eqnhtt:3b}
		\end{numcases}
	\end{subequations}
\end{small}
where
\begin{eqnarray*}
	&&a_{1}(\H,\F)=\int_{\Omega}{\d{\ep}_{r}^{-1}\curl\H\cdot\curl\F^{*}}\mathrm{d}\Omega \\
	&&b_{1}(\F,q)=\int_{\Omega}{\nabla q\cdot\F^{*}}\mathrm{d}\Omega\\
	&&c_{1}(\H,\upsilon)=\int_{\Omega}{\d{\mu}_{r}\H\cdot{\grad{\upsilon^{*}}}}\mathrm{d}\Omega\\
	&&d_{1}(\H,\F)=\int_{\Omega}{\d{\mu}_{r}\H\cdot\F^{*}}\mathrm{d}\Omega.
\end{eqnarray*}
The mixed variational form associated with PDEs (\ref{eqq:2}) reads as:
\begin{small}
\begin{subequations} \label{eqntt:3}
	\begin{numcases}{}
		\textrm{Seek~} \Lambda\in{\mathbb{C}},~{\bf{0}}\neq\E\in{\mathbb{H}_{0}(\textrm{curl},\Omega)},~p\in{H_{0}^1(\Omega)},\textrm{~such that~}\nonumber\\
		a_{2}(\E,\F)+\alpha b_{2}(\F,p)= \Lambda d_{2}(\E,\F),~\forall\F\in{\mathbb{H}_{0}(\textrm{curl},\Omega)}\label{ttformm1}\\
		c_{2}(\E,\upsilon) = 0,~\forall \upsilon\in{H_{0}^1(\Omega)}\label{ttformm2}
	\end{numcases}
\end{subequations}
\end{small}
where
\begin{eqnarray*}
	&&a_{2}(\E,\F)=\int_{\Omega}{\d{\mu}_{r}^{-1}\curl\E\cdot\curl\F^{*}}\mathrm{d}\Omega \\
	&&b_{2}(\F,p)=\int_{\Omega}{\nabla p\cdot\F^{*}}\mathrm{d}\Omega\\
	&&c_{2}(\E,\upsilon)=\int_{\Omega}{\d{\ep}_{r}\E\cdot{\grad{\upsilon^{*}}}}\mathrm{d}\Omega\\
	&&d_{2}(\E,\F)=\int_{\Omega}{\d{\ep}_{r}\E\cdot\F^{*}}\mathrm{d}\Omega.	
\end{eqnarray*}

Let $\mathcal{T}_{h}$ be a regular tetrahedral mesh in $\Omega$, where $h$ is the length of the longest edge
in $\mathcal{T}_{h}$. For the sake of simplicity, here the edge space $\W_{1}^{h}$ of the lowest order under the mesh $\mathcal{T}_{h}$ is used to approximate the Hilbert space $\mathbb{H}(\textrm{curl},\Omega)$. From \cite{monk2003finite}, the vector function space $\W_{1}^{h}$ is a finite element subspace of $\mathbb{H}(\textrm{curl},\Omega)$ and it implies that $\W_{1}^{h}$ is a curl-conforming space. Set $\W_{2}^{h}=\W_{1}^{h}\bigcap\mathbb{H}_{0}(\textrm{curl},\Omega)$,
then it is known that $\W_{2}^{h}$ is a good approximation of $\mathbb{H}_{0}(\textrm{curl},\Omega)$. In addition, the standard linear element space $S_{1}^{h}$ is applied to approximate the Hilbert space $H^{1}(\Omega)$. Set $S_{2}^{h}=S_{1}^{h}\bigcap H_{0}^{1}(\Omega)$, then it is known that $S_{2}^{h}$ is a good approximation of $H_{0}^{1}(\Omega)$. The definition of $\W_{1}^{h}$ is given in \cite{Jiang2016}. For details, please see \cite{Jiang2016}.

Restricting the mixed variational form (\ref{eqnhtt:3}) on  $\W_{1}^{h}\times S_{1}^{h}$, one can get
the discrete mixed variational form associated with (\ref{eqnhtt:3}), and that is
\begin{small}
	\begin{subequations} \label{eqnH}
		\begin{numcases}{}
			\textrm{Seek~} \Lambda_{h}\in{\mathbb{C}},~{\bf{0}}\neq\H_{h}\in{\W_{1}^{h}},
			~q_{h}\in{S_{1}^{h}} \textrm{~such that}\nonumber\\
			a_{1}(\H_{h},\F)+\beta b_{1}(\F,q_{h})= \Lambda_{h} d_{1}(\H_{h},\F),\forall\F\in{\W_{1}^{h}}\label{eqnH:1}\\
			c_{1}(\H_{h},\upsilon) = 0,~~\forall\upsilon\in{S_{1}^{h}}\label{eqnH:2}
		\end{numcases}
	\end{subequations}
\end{small}

Restricting the mixed variational form (\ref{eqntt:3}) on $\W_{2}^{h}\times S_{2}^{h}$, one can get
the discrete mixed variational form associated with (\ref{eqntt:3}), and that is
\begin{small}
	\begin{subequations} \label{eqnE}
		\begin{numcases}{}
			\textrm{Seek~} \Lambda_{h}\in{\mathbb{C}},~{\bf{0}}\neq\E_{h}\in{\W_{2}^{h}},
			~p_{h}\in{S_{2}^{h}} \textrm{~such that}\nonumber\\
			a_{2}(\E_{h},\F)+\alpha b_{2}(\F,p_{h})= \Lambda_{h} d_{2}(\E_{h},\F),\forall\F\in{\W_{2}^{h}}\label{eqnE:1}\\
			c_{2}(\E_{h},\upsilon) = 0,~~\forall\upsilon\in{S_{2}^{h}}\label{eqnE:2}
		\end{numcases}
	\end{subequations}
\end{small}

\begin{table*}[ht]
\renewcommand{\arraystretch}{1.3}
\centering
\caption{\label{ser1} The Physical Eigenvalue ($\mathrm{m}^{-2}$) from A Cylindrical Resonant Cavity}
\begin{tabular}{ccccc}
 \hline
 $h(\textrm{m})$& 0.1043  & 0.0714 & 0.0580 &0.0428 \\
 \hline
  $\Lambda_{h}(\textrm{cu},\E)$&  $24.5133 - 7.5592\textrm{j}$ &$24.3328 - 7.5559\textrm{j}$ & $24.2950- 7.5599\textrm{j}$  &$24.2499-7.5594\textrm{j}$\\
   $t_{1}(s)$& 9  & 60&   110&490\\
 \hline
 $\Lambda_{h}(\textrm{cu},\H)$&  $24.5129 - 7.5588\textrm{j}$ &$24.3321 - 7.5550\textrm{j}$ & $24.2944- 7.5591\textrm{j}$  &$24.2490-7.5585\textrm{j}$\\
 $t_{2}(s)$& 10  & 65&   128&550\\
 \hline
  $\Lambda_{h}(\textrm{pr},\H)$~\cite{Jiang2020}&  $24.5131 - 7.5590\textrm{j}$ &$24.3324 - 7.5554\textrm{j}$ & $24.2947 - 7.5595\textrm{j}$  &$24.2495-7.5589\textrm{j}$\\
   $t_{3}(s)$&30   &180 &  1002 &12000\\
   \hline
\end{tabular}
\end{table*}

\begin{table*}[ht]
\renewcommand{\arraystretch}{1.3}
 \centering
 \caption{\label{TEWscalar} The Physical Eigenvalue ($\mathrm{m}^{-2}$) from A Torus Resonant Cavity}
 \begin{tabular}{ccccc}
 \hline
 $h(\mbox{m})$
 &  0.3465 & 0.2680 & 0.1882& 0.1322\\
\hline
$\Lambda_{h}(\textrm{cu},\E)$&     $1.0938 + 3.6441\textrm{j}$ &$1.0720+ 3.5728\textrm{j}$ & $1.0509+3.5331\textrm{j}$  &$1.0396+ 3.5060\textrm{j}$\\
$\Lambda_{h}(\textrm{cu},\H)$&  $1.0932 + 3.6435\textrm{j}$ &$1.0712+ 3.5720\textrm{j}$ & $1.0499+3.5320\textrm{j}$  &$1.0382+ 3.5049\textrm{j}$\\
$\Lambda_{h}(\textrm{COMSOL},\E)$&  $1.0936 + 3.6438\textrm{j}$ &$1.0716+ 3.5724\textrm{j}$ & $1.0505+3.5326\textrm{j}$  &$1.0389+ 3.5053\textrm{j}$\\
\hline
\end{tabular}
\end{table*}

$(\Lambda_{h},\H_{h},q_{h})$ in (\ref{eqnH}) and $(\Lambda_{h},\E_{h},p_{h})$ in (\ref{eqnE}) are an approximations of the exact solution $(\Lambda,\H,q)$ in (\ref{eqnhtt:3}) and $(\Lambda,\E,p)$ in (\ref{eqntt:3}) respectively. Next, the MFEM is applied to solve the discrete mixed variational forms (\ref{eqnH}) and (\ref{eqnE}).
\subsection{Generalized Eigenvalue Problem}
Suppose that $S_{1}^{h}=\textrm{span}\{s_{i}\}_{i=1}^{n}$ and $\W_{1}^{h}=\textrm{span}\{\mathcal{W}_{i}\}_{i=1}^{m}$, where $s_{i}$ and $\mathcal{W}_{i}$ is $i$th
the global basis functions of $S_{1}^{h}$ and $\W_{1}^{h}$. The local basis functions of the finite element spaces $S_{1}^{h}$  and $\W_{1}^{h}$ is given in \cite{Jiang2016}.
It's important to point out these three boundary conditions (\ref{eqq:3c})-(\ref{eqq:3e}) are three natural boundary conditions in FEM. Hence, the boundary conditions (\ref{eqq:3c})-(\ref{eqq:3e}) do not need a special treatment in the numerical implementation of MFEM.
However, the boundary conditions (\ref{eqq:2c})-(\ref{eqq:2d}) need a special treatment in the numerical implementation of MFEM, since the boundary conditions (\ref{eqq:2c})-(\ref{eqq:2d}) are two essential boundary condition in FEM. Here, we only list the matrix eigenvalue problem corresponding to the mixed variational form (\ref{eqnH}), and that is
\begin{equation}\label{eig1}
    \left[
      \begin{array}{cc}
        A & \beta B \\
        C & O \\
      \end{array}
    \right]\left[
             \begin{array}{c}
               \xi \\
               \zeta \\
             \end{array}
           \right]
    =\Lambda_{h}\left[
      \begin{array}{cc}
        D & O \\
        O & O \\
      \end{array}
    \right]\left[
             \begin{array}{c}
               \xi \\
               \zeta \\
             \end{array}
           \right],
\end{equation}
where
\begin{gather*}
   \xi=[\xi_{1},\xi_{2},\cdots,\xi_{n}]^{T}, \quad \zeta=[\zeta_{1},\zeta_{2},\cdots,\zeta_{m}]^{T}\\
   A=(a_{ik})\in{\mathbb{C}^{n\times n}}, \quad
   C=(c_{ik})\in{\mathbb{C}^{m\times n}}\\
   B=(b_{ik})\in{\mathbb{C}^{n\times m}}, \quad
   D=(d_{ik})\in{\mathbb{C}^{n\times n}}\\
   a_{ik}=a_{1}(\mathcal{W}_{k},\mathcal{W}_{i}), \quad
   c_{ik}=c_{1}(\mathcal{W}_{k},s_{i})\\
   b_{ik}=b_{1}(\mathcal{W}_{k},s_{i}), \quad
   d_{ik}=d_{1}(\mathcal{W}_{k},\mathcal{W}_{i})
\end{gather*}

According to actual calculations, we find that the infinite norm from the matrices $A$ and $B$ are quite different.
Based on this reason, we take $\beta=\varrho(A)/\varrho(B)$ in (\ref{eig1}), where $\varrho(A)$ and $\varrho(B)$ are the infinite norm of the matrices $A$ and $B$, respectively. In (\ref{eqnE}), the constant $\alpha$ is also selected in this way.

After solving the generalized matrix eigenvalue problem (\ref{eig1}) by the implicitly restarted Arnoldi methods in package ARPACK \cite{Lehoucq}, we can obtain several numerical physical eigenvalues.

\section{Numerical experiments}
In this section, two numerical experiments are carried out to show that the numerical method given in the article can remove all the spurious modes, including spurious dc modes.

Because it is very difficult to solve the resonant cavity problem filled with anisotropic media by an explicit analytical method, numerical methods must be used to solve this problem. How to measure whether the numerical result of a numerical solution is correct or not? One way is that we can select the numerical examples in the existing references as benchmarks. The other way is that we can simulate this problem by a commercial software, such as COMSOL Multiphysics, ANSYS and so on.
Here, we employ COMSOL Multiphysics to simulate the 3-D closed resonant cavity problems and then compare with the numerical results from COMSOL Multiphysics and the method given in the article. The computational model in COMSOL Multiphysics is as follows:
\begin{subequations} \label{eqcom:2}
\begin{numcases}{}
  \curl\Big({\d{\mu}_{r}^{-1}}\curl\E\Big) =\Lambda{\d{\ep}}_{r}\E~~\text{in}~\Omega\label{eqcom:2a}\\
  \^n\times\E = {\bf{0}}~~\text{on}~\p\Omega\label{eqcom:2c}
\end{numcases}
\end{subequations}
It is clear that PDEs (\ref{eqcom:2}) ignores the divergence-free condition (\ref{eq:2b}). As a consequence, there are many spurious dc eigenvalues in the numerical eigenvalues from COMSOL Multiphysics. In view of the fact that the spurious dc eigenvalues are well distinguished from the numerical eigenvalues of COMSOL Multiphysics, we can select physical eigenvalues in the numerical eigenvalues of COMSOL Multiphysics.

Our MFEM codes are written in MATLAB, and they run on a computer with 8-cores Intel Core i9-11900K 3.5GHz CPU and 128GB-RAM.
The mesh data of our MFEM codes are from the ones of COMSOL. In the numerical simulations of COMSOL, we select the FEM based on the edge element of the lowest order to simulate PDEs (\ref{eqcom:2}).

In order to distinguish the numerical eigenvalues from (\ref{eqnH}) and (\ref{eqnE}). The numerical eigenvalues from (\ref{eqnH}) and (\ref{eqnE}) are denoted by $\Lambda_{h}(\textrm{cu},\H)$ and $\Lambda_{h}(\textrm{cu},\E)$, respectively.
The numerical eigenvalues obtained by the projection method in \cite{Jiang2020} to solve PDEs (\ref{eq:3}) are denoted by
$\Lambda_{h}(\textrm{pr},\H)$. In addition, the numerical eigenvalues obtained by COMSOL Multiphysics to solve PDEs (\ref{eqcom:2}) are denoted by $\Lambda_{h}(\textrm{COMSOL},\E)$.

\subsection{Cylindrical Cavity}
The radius and height of cylindrical cavity is denoted by $r$ and $b$, respectively. Set $r=0.2$\,m and $b=0.5$\,m. Suppose  that the relative permittivity and permeability tensor of the anisotropic medium in the whole cylindrical cavity are
\begin{equation*}
    \d{\ep}_{r}=
    \begin{bmatrix}
    2+\textrm{j}&0&0\\
    0&2+\textrm{j}&0\\
    0&0&2
    \end{bmatrix},\quad
    \d{\mu}_{r}=
    \begin{bmatrix}
        2-\textrm{j}&0.375\textrm{j}&0\\
    0.375\textrm{j}&2-\textrm{j}&0\\
    0&0&2
    \end{bmatrix}.
\end{equation*}
It is clear that both $\d{\ep}_{r}^{\dag}\neq\d{\ep}_{r}$ and $\d{\mu}_{r}^{\dag}\neq\d{\mu}_{r}$ are valid, and this implies that the anisotropic medium is both electric and magnetic lossy.

In \cite{Jiang2020}, the penalty method, augmented method and projection method are used to solve the above numerical example and the solved model is just PDEs (\ref{eq:3}). Among these three computational methods, only the projection method does not introduce any spurious modes.
Hence, in Table \ref{ser1}, we only list the numerical eigenvalue $\Lambda_{h}(\textrm{pr},\H)$ of the projection method.
The time $t_{1}$  and $t_{2}$ are obtained by $\Lambda_{h}(\textrm{cu},\E)$ and $\Lambda_{h}(\textrm{cu},\H)$ using the computational methods in the current article respectively. The time $t_{3}$ are obtained  by the projection method in \cite{Jiang2020}. For simplicity, we only list the numerical eigenvalue corresponding to dominate mode in Table \ref{ser1}.

As shown in Table \ref{ser1}, $\Lambda_{h}(\textrm{cu},\E)$ and $\Lambda_{h}(\textrm{cu},\H)$ are approximately equal to $\Lambda_{h}(\textrm{pr},\H)$.
$t_{1}\approx t_{2}<t_{3}$ implies that the computational method supported by this article is more efficient than the projection method in \cite{Jiang2020}. The reason for the low computational efficiency of the projection method is that this projection method uses a complete singular value decomposition algorithm to find the dense basis vector of the null space of a rectangular matrix. In contrast, the algorithm in this article does not involve the calculation of dense matrices. It is worth pointing out that the algorithm provided by this article can remove all the spurious modes when it is used to solve the resonant cavity problem filled with a both electric and magnetic lossy, anisotropic medium.

\subsection{Torus Cavity}
The parametric definition of a torus is as follows:
\begin{subequations}
\begin{numcases}{}
x(\theta,\varphi) =(\rho_{1}+\rho_{2}\cos\theta)\cos\varphi,~~~~\rho_{1}>\rho_{2}>0\nonumber\\
y(\theta,\varphi)=(\rho_{1}+\rho_{2}\cos\theta)\sin\varphi,~~~~\theta,~\varphi\in[0,2\pi)\nonumber\\
z(\theta,\varphi)=\rho_{2}\sin\theta,\nonumber
\end{numcases}
\end{subequations}
where $\rho_{1}$ and $\rho_{2}$ are major and minor radiuses, respectively. Set $\rho_{1}=0.8$\,m, $\rho_{2}=0.4$\,m. Assume that the torus cavity is filled with a homogeneous, both electric and magnetic lossy, anisotropic medium. The medium parameter of the material in this torus cavity is given by
\begin{equation*}
    \d{\ep}_{r}=
    \begin{bmatrix}
    2-0.5\textrm{j}&0.25\textrm{j}&0.25\textrm{j}\\
    -0.25\textrm{j}&2-0.5\textrm{j}&0.25\textrm{j}\\
    -0.25\textrm{j}&-0.25\textrm{j}&2-0.5\textrm{j}
    \end{bmatrix}
\end{equation*}

\begin{equation*}
    \d{\mu}_{r}=
    \begin{bmatrix}
        1-0.2\textrm{j}&0&0\\
    0&1-0.4\textrm{j}&0\\
    0&0&1-0.8\textrm{j}
    \end{bmatrix}.
\end{equation*}
This torus cavity problem is simulated by the numerical methods in this article, the projection method in \cite{Jiang2020} and COMSOL Multiphysics.

In Table \ref{TEWscalar}, we list the numerical eigenvalues $\Lambda_{h}(\textrm{cu},\E)$, $\Lambda_{h}(\textrm{cu},\H)$ and $\Lambda_{h}(\textrm{COMSOL},\E)$ associated with the dominant mode. We can find that $\Lambda_{h}(\textrm{cu},\E)\approx\Lambda_{h}(\textrm{cu},\H)\approx\Lambda_{h}(\textrm{COMSOL},\E)$. This implies that the implementation of our MATLAB code is correct. It is worthwhile to point out that there are not any numerical spurious  results from our MATLAB codes, however, there are many spurious dc eigenvalues from COMSOL.

In addition, according to the above two numerical examples $A$ and $B$, we do find that the entries from $\zeta$ in (\ref{eig1}) are all the same and $p_{h}\approx0$ in (\ref{eqnE:1}) is valid. These verify the correctness of the theory in Section \uppercase\expandafter{\romannumeral2}.

\section{Conclusion}
By introducing dummy variables, this article enforces the divergence-free condition supported by Gauss's law in electromagnetics in a weak sense. As a consequence, this novel method can be free of all the spurious modes in solving the resonant cavity problem filled with both electric and magnetic lossy, anisotropic media.

In the current article, the MFEM based upon the edge element of the lowest order and the standard linear element is to find the eigenmodes of the resonant cavity. It is known that the convergence rate of this type of MFEM is very slow, therefore, we need to accelerate the convergence rate of the MFEM. This leads to the MFEM based on higher order basis functions.
In the future, the MFEM based upon higher order basis functions is applied to solve 3-D classic resonant cavity problem.
Notice that this MFEM also has an ability of removing all the spurious modes.

\end{document}

%% file: latex_defs.tex
\def\^{\hat}
\def\~{\tilde}

\def\3h{{3\over 2}}

\def\v#1{{\bf#1}}

\def\E{{\v E}}

\def\W{{\v W}}

\def\H{{\v H}}

\def\F{{\v F}}

\def\ep{\epsilon}

\def\p{\partial}
\def\pd#1#2{{{\partial #1}\over{\partial #2}}}

\def\div{\nabla\cdot}

\def\curl{\nabla\times}
\def\grad{\nabla}

\def\eqn#1$${\eqno{{\rm #1}}$$}

\def\ep{\epsilon}

\def\~{\tilde}

\def\^{\hat}

\def\d#1{\overline{\overline{#1}}}

\def\XXint#1#2#3{{\setbox0=\hbox{$#1{#2#3}{\int}$}
     \vcenter{\hbox{$#2#3$}}\kern-.5\wd0}}
